\documentclass[10pt,twoside]{article}
\usepackage{graphicx}
\usepackage{amsmath}
\usepackage{amssymb}
\usepackage{Latex-document}

\def\volumeno{I}

\def\underTilde#1{{\baselineskip=0pt\vtop{\hbox{$#1$}\hbox{$\sim$}}}{}}

\def\({\left(}
\def\){\right)}

\def\{{\left\lbrace}
\def\}{\right\rbrace}

\def\<{\langle}
\def\>{\rangle}

\def\qedbox{$\rlap{$\sqcap$}\sqcup$}
\def\qed{\hbox to 0pt{}\nobreak\hfill\qedbox}

\newcommand{\sqed}{
  \ifmmode\eqno{\hbox{\qedbox}}
  \else\hbox to 0pt{}\nobreak\hfill\qedbox
  \fi
  }

\begin{document}


\markboth{Beyond $\underTilde{\Sigma}^2_1$ Absoluteness}{W. Hugh Woodin}
\title{\bf \boldmath Beyond $\underTilde{\Sigma}^2_1$ Absoluteness\vskip 6mm}
\author{W. Hugh Woodin\vspace*{-0.5cm}\thanks{Department of Mathematics,
University of California, Berkeley, Berkeley CA 94720, USA. E-mail:
woodin@math.berkeley.edu}}
\date{\vspace{-8mm}}

\maketitle

\thispagestyle{first} \setcounter{page}{515}

\begin{abstract}

\vskip 3mm

There have been many generalizations of Shoenfield's Theorem
on the absoluteness of $\Sigma^1_2$ sentences between uncountable
transitive models of $\mathrm{ZFC}$. One of the strongest versions
currently known deals with $\Sigma^2_1$ absoluteness conditioned
on $\mathrm{CH}$. For a variety of reasons, from the study of inner
models and from simply
combinatorial set theory, the question of
whether conditional $\Sigma^2_2$ absoluteness is possible at all,
and if so, what large cardinal assumptions are involved and what
sentence(s) might play the role of $\mathrm{CH}$, are
fundamental questions.  This article investigates the possiblities for
$\Sigma^2_2$ absoluteness by extending the connections between
determinacy hypotheses and absoluteness hypotheses.

\vskip 4.5mm


\end{abstract}

\vskip 12mm

\section{Absoluteness and strong logics} \label{section 1}\setzero
\vskip-5mm \hspace{5mm }

There have been many generalizations of Shoenfield's Theorem
on the absoluteness of $\Sigma^1_2$ sentences between uncountable
transitive models of $\mathrm{ZFC}$. Absoluteness theorems are
meta-mathematically interesting since they identify levels of complexity
where the technique of forcing cannot be used to establish independence.

A sentence, $\phi$, is a
{\em $\Sigma^2_1$-sentence} if for some $\Sigma_1$-formula, $\psi(x)$,
$\phi$ is provably equivalent in $\mathrm{ZFC}$, Zermelo Frankel set theory with
the {\em Axiom of Choice},  to the
assertion that $\psi[\mathbb{R}]$ holds. While this is not the standard
definition, for the purposes of this article it is equivalent.

\medskip
{\bf Theorem 1.1}  {\em Suppose that $\phi$ is a $\Sigma^2_1$
sentence, there exists a proper class of measurable Woodin
cardinals and that $\mathrm{CH}$ holds. Suppose $\mathbb{P}$ is a
partial order and that $V^{\mathbb P} \vDash \mathrm{CH}$. Then $V
\vDash \phi$ if and only if $V^{\mathbb{P}} \vDash \phi $}\qed
\medskip

This is $\Sigma^2_1$ generic absoluteness conditioned on $\mathrm{CH}$. Because
$\mathrm{CH}$ is itself a $\Sigma^2_1$ sentence, this {\em conditional}
form of $\Sigma^2_1$ generic absoluteness is the best one can hope for.
The meta-mathematical significance of this kind of absoluteness result is simply this.
If a problem is expressible as a $\Sigma^2_1$ sentence, and there are many such examples
from analysis, then it is likely settled by $\mathrm{CH}$ (augmented by modest large
cardinal hypotheses).  The technique of forcing
cannot be used to demonstrate otherwise.

Absoluteness theorems can be naturally reformulated
using {\em strong logics}. For generic absoluteness the
relevant logic is $\Omega^*$-logic.

\medskip
{\bf Definition 1.2}{\bf ($\Omega^*$-logic)} {\it Suppose that there
exists a proper class of Woodin cardinals and that $\phi$ is a sentence.
Then
$$
\mathrm{ZFC} \vdash_{\Omega^*} \phi
$$
if for all ordinals $\alpha$ and for all partial orders $\mathbb{P}$ if
$V_{\alpha}^{\mathbb{P}} \vDash \mathrm{ZFC}$, then $
V_{\alpha}^{\mathbb{P}}\vDash \phi$}.\qed
\medskip

The theorem on $\Sigma^2_1$-absoluteness and $\mathrm{CH}$ can be
reformulated as follows.

\medskip
{\bf Theorem 1.3} {\em Suppose there exists a proper class of
measurable Woodin cardinals. Then for each $\Sigma^2_1$
sentence $\phi$, either
$\mathrm{ZFC} + \mathrm{CH} \vdash_{\Omega^*} \phi$;
or
$\mathrm{ZFC} + \mathrm{CH} \vdash_{\Omega^*} (\neg \phi)$}.\qed
\medskip

But there is another natural strong logic; $\Omega$-logic,
the definition of $\Omega$-logic involves
universally Baire sets of reals.

\medskip
{\bf Definition 1.4}{ \cite{fmw}} {\it A set $A \subseteq \mathbb{R}^n$
is {\em universally Baire} if for any continuous function, $F:\Omega \to
\mathbb{R}^n$, where $\Omega$ is a compact Hausdorff space, the preimage
of $A$,
$$
\{p \in X{\large\;\mid\;} F(p) \in A\},
$$
has the property of Baire in $\Omega$; \mbox{i.\ e.} is open in $\Omega$
modulo a meager set}.\qed
\medskip

Every borel set $A \subseteq \mathbb{R}^n$
is universally Baire. More generally
the universally Baire sets form a $\sigma$-algebra
closed under preimages by borel functions
$$
f:\mathbb{R}^n \to \mathbb{R}^m.
$$
The universally Baire sets have the classical regularity
properties of the borel sets, for example they are Lebesgue measurable
and have the property of Baire.
If there exists a proper class of Woodin cardinals then
the universally Baire sets are closed under projection
and every universally Baire set is determined.

Suppose that
$A \subseteq \mathbb{R}$
in universally Baire and that $V[G]$ is a set generic extension
of $V$. Then the set $A$ has canonical interpretation
as a set
$$
A_G \subseteq \mathbb{R}^{V[G]}.
$$
The set $A_G$ is defined as
$$
A_G = \cup\{\mbox{\rm range}(\pi_G){\large\;\mid\;}
\pi \in V, \mbox{\rm range}(\pi) = A\};
$$
here $\pi$ is a function, $\pi:\lambda^{\omega} \to \mathbb{R}$,
which satisfies the uniform continuity requirement that
for $f \neq g$;
$$
|\pi(f) - \pi(g)| < 1/(n+1)
$$
where $n <\omega$ is least such that $f(n) \neq g(n)$.
 If there exists a proper class of Woodin cardinals
then
$$
\<H(\omega_1),A,\in\> \prec \<H(\omega_1)^{V[G]},A_G,\in\>.
$$

\medskip
{\bf Definition 1.5} {\it Suppose that $A \subseteq \mathbb{R}$ is
universally Baire and that $M$ is a transitive set such that $M \vDash
\mathrm{ZFC}$. Then $M$ is {\em $A$-closed} if for each partial order
$\mathbb{P} \in M$, if $G \subset \mathbb{P}$ is $V$-generic then in
$V[G]$: $A_G \cap M[G] \in M[G]$}.\qed
\medskip

\medskip
{\bf Definition 1.6}{\bf ($\Omega$ logic)} {\it Suppose that there
exists a proper class of Woodin cardinals and that $\phi$ is a sentence.
Then $\mathrm{ZFC} \vdash_{\Omega} \phi$ if there exists a universally
Baire set $A \subseteq \mathbb{R}$ such that if $M$ is any countable
transitive set such that $M \vDash \mathrm{ZFC}$ and such that $M$ is
$A$-closed, then $ M\vDash \phi$}.\qed
\medskip

Both $\Omega$-logic and $\Omega^*$-logic are definable and
generically invariant.

A natural question, given the theorem on
generic absoluteness for $\Sigma^2_1$ is the following
question:

\begin{quote} {\em Suppose there exists a proper class of
measurable Woodin cardinals. Does it
follow that for each $\Sigma^2_1$
sentence $\phi$, either
$\mathrm{ZFC} + \mathrm{CH} \vdash_{\Omega} \phi;$
or
$\mathrm{ZFC} + \mathrm{CH} \vdash_{\Omega} (\neg \phi)?$}
\end{quote}

The answer is yes if ``iterable'' models with
measurable Woodin cardinals exist.

\medskip

{\bf $\Omega$ Conjecture:}
\begin{quote}{\em
Suppose that there exists a proper class of Woodin
cardinals and that $\phi$ is a $\Pi_2$ sentence. Then
$\mathrm{ZFC} \vdash_{\Omega^*} \phi$
if and only if $\mathrm{ZFC} \vdash_{\Omega} \phi$.}
\end{quote}

\medskip

It is immediate
from the definitions and Theorem~1.1, that the
$\Omega$~Conjecture settles the question above affimatively.
But the consequences of the $\Omega$-Conjecture are far more
reaching. If the $\Omega$~Conjecture is true, then generic absoluteness is
equivalent to absoluteness in $\Omega$-logic and this in turn has
significant metamathematical implications.

We fix some conventions.
A formula, $\phi(x)$, is  a
{\em $\Sigma^2_2$-formula} if for some $\Sigma_2$-formula, $\psi(x)$,
the formula $\phi(x)$ is provably equivalent in $\mathrm{ZFC}$ to the
formula
$$
`` x \in H(c^+) \mbox { and } \<H(c^+), \in\> \vDash \psi[x]".
$$

Finally  $\phi(x)$ is  a
{\em $\Sigma^2_2({\cal I}_{_{\mathrm{NS}}})$-formula} if for some $\Sigma_2$-formula, $\psi(x)$,
the formula $\phi(x)$ is provably equivalent in $\mathrm{ZFC}$ to the
formula
$$
`` x \in H(c^+) \mbox { and } \<H(c^+),{\cal
I}_{_{\mathrm{NS}}},\in\> \vDash \psi[x]".
$$
where ${\cal I}_{_{\mathrm{NS}}}$ denotes the nonstationary ideal on $\omega_1$.

There is a limit to the possible extent of absoluteness in $\Omega$-logic.
One version is given by the following theorem.

\medskip
{\bf Theorem 1.7} {\em Suppose that there exist
a proper class of Woodin cardinals, $\Psi$ is a sentence
and that for each $\Sigma^2_2({\cal I}_{_{\mathrm{NS}}})$ sentence $\phi$, either
$\mathrm{ZFC} + \Psi
\vdash_{\Omega} \phi$, or
$\mathrm{ZFC} + \Psi
\vdash_{\Omega} (\neg\phi)$.
Then ${\mathrm{ZFC}} + \Psi$ is $\Omega$-inconsistent. \qed}
\medskip

In short:
\begin{quote}{\em $\Sigma^2_2({\cal I}_{_{\mathrm{NS}}})$ absoluteness is not possible in
$\Omega$-logic. If the $\Omega$~Conjecture holds then generic
absoluteness for $\Sigma^2_2({\cal I}_{_{\mathrm{NS}}})$ sentences is not possible.}
\end{quote}

So for absoluteness in $\Omega$-logic the most one can hope for is
that there exist a sentence $\Psi$ such that for each $\Sigma^2_2$
sentence $\phi$, either $\mathrm{ZFC} + \Psi \vdash_{\Omega}
\phi$, or
 $\mathrm{ZFC} + \Psi
\vdash_{\Omega} (\neg\phi)$.
In particular if the $\Omega$~Conjecture holds then $\Sigma^2_2$ generic absoluteness
is the most one can hope for.

Suppose that $\Psi$ is a sentence such
that for each $\Sigma^2_2$ sentence $\phi$, either
$\mathrm{ZFC} + \Psi
\vdash_{\Omega} \phi$, or
$\mathrm{ZFC} + \Psi
\vdash_{\Omega} (\neg\phi)$.
Then
$\mathrm{ZFC} + \Psi \vdash_{\Omega} 2^{\aleph_0} < 2^{\aleph_1}$.
A natural conjecture is that in fact
$\mathrm{ZFC} + \Psi \vdash_{\Omega} \mathrm{CH}$.

In any case from this point on we shall consider absoluteness
just in the context of $\mathrm{CH}$.

Generic absoluteness is closely related to determinacy. The
statement of  a theorem which illustrates one aspect of
this requires the following definition.

\medskip
{\bf Definition 1.8} {\it Suppose that there exists a proper class of
Woodin cardinals. A set $A \subseteq \mathbb{R}$ is $\Omega^*$-{\em
recursive} if there exists a formula $\phi(x)$ such that:
\begin{enumerate}

\item $A = \{r {\large\;\mid\;} \mathrm{ZFC} \vdash_{\Omega^*}\phi[r]\}$;

\item For all
partial orders, $\mathbb{P}$, if $G \subseteq \mathbb{P}$
is $V$-generic then for each $r\in \mathbb{R}^{V[G]}$,
either
$$
V[G] \vDash \mbox{\rm ``$\mathrm{ZFC} \vdash_{\Omega^*} \phi[r]$''},
$$
or
$V[G] \vDash \mbox{\rm ``$\mathrm{ZFC} \vdash_{\Omega^*} (\neg\phi)[r]$''}$.\qed
\end{enumerate}}
\medskip

{\bf Theorem 1.9} {\em Suppose that there exists
a proper class of Woodin cardinals.  Suppose
that $A \subseteq \mathbb{R}$ is $\Omega^*$-recursive.
Then $A$ is determined.\qed}
\medskip

On the other hand there are many examples where suitable
determinacy assumptions imply generic absoluteness.
Our main results deal with generalizations of these connections
to $\Sigma^2_1$ and $\Sigma^2_2$ in the context of $\mathrm{CH}$.

\section{Absoluteness and determinacy} \label{section 2}\setzero
\vskip-5mm \hspace{5mm }

We fix a reasonable coding of elements of $H(\omega_1)$ by reals.
This is simply a surjection
$$
\pi:\mathrm{dom}(\pi) \to H(\omega_1)
$$
where $\mathrm{dom}(\pi) \subseteq \mathbb{R}$.
All we require of $\pi$ is that $\pi\in L(\mathbb{R})$; the
natural choice for $\pi$ is definable in $H(\omega_1)$.
For each set $X \subseteq H(\omega_1)$ let
$$
X^* = \{x \in \mathbb{R}{\large\;\mid\;} \pi(x) \in X\}.
$$

Suppose $X \subseteq \{0,1\}^{\omega_1}$. Associated to $X$
is a game of length $\omega_1$.  The convention is that Player~I
plays first at limit stages.  Strategies are functions:
$$
\tau: \{0,1\}^{< \omega_1} \to \{0,1\}.
$$
Suppose that $\Gamma \subseteq {\cal P}(\mathbb{R})$. Then $X$ is {\em $\Gamma$-clopen} if
there exist sets $Y \subset H(\omega_1) $  and $Z \subset H(\omega_1)$
such that
\begin{enumerate}
\item $Y\cap Z = \emptyset$,
\item for all $a \in \{0,1\}^{\omega_1}$ there exists $\alpha < \omega_1$
such that either $a|\alpha \in Y$ or $a|\alpha \in Z$,
\item $X$ is the set of $a \in \{0,1\}^{\omega_1}$ such that there exists
$\alpha < \omega_1$ such that $a|\alpha \in Y$ and such that
$a|\beta \notin Z$ for all $\beta < \alpha$,
\item $Y^* \in \Gamma$ and $Z^* \in \Gamma$.
\end{enumerate}

The first result on the determinacy of $\Gamma$-clopen sets
is due to Itay Neeman. One version is the following.

\medskip
{\bf Theorem 2.1} \cite{neemanb} {\em Suppose that there exists a Woodin
cardinal which is a limit of Woodin cardinals. Suppose that
$X \subseteq \{0,1\}^{\omega_1}$ and $X$ is $\Pi^1_1$-clopen.

Then $X$ is determined.\qed}
\medskip

The proof of Neeman's theorem combined with techniques
from the fine structure theory associated to $\mathrm{AD}^+$
yields the following generalization which we shall need.

\medskip
{\bf Theorem 2.2} {\em Suppose that there is a proper class
of  Woodin cardinals which are limits of Woodin cardinals.  Let
$\Gamma^{\infty}$ be the set of all $A \subseteq \mathbb{R}$
such that $A$ is universally Baire and suppose that
$X \subseteq \{0,1\}^{\omega_1}$ is such that $X$ is $\Gamma^{\infty}$-clopen.

Then $X$ is determined.\qed}
\medskip

Suppose $X \subseteq \{0,1\}^{\omega_1}$
and that $\Gamma \subseteq {\cal P}(\mathbb{R})$.
Then $X$ is {\em $\Gamma$-open} if
there exist sets $Y \subset H(\omega_1) $
such that
\begin{enumerate}
\item $X$ is the set of $a \in \{0,1\}^{\omega_1}$ such that there exists
$\alpha < \omega_1$ such that $a|\alpha \in Y$ .
\item $Y^* \in \Gamma$.
\end{enumerate}

John Steel has proved that under fairly general conditions, if
$\Gamma \subseteq {\cal P}(\mathbb{R})$ is such that all $\Gamma$-open
sets are determined then for each $X \subseteq \{0,1\}^{\omega_1}$,
if $X$ is $\Gamma$-open and if Player~I wins the game given
by $X$, then there is a winning strategy for Player~I which is
(suitably) definable from parameters in $\Gamma$; \cite{steel}. The following is a straightforward
corollary:

\medskip
{\bf Corollary 2.3} {\em Suppose that there exists a proper class
 of Woodin cardinals.
Let
$\Gamma^{\infty}$ be the set of all $A \subseteq \mathbb{R}$
such that $A$ is universally Baire and suppose that for each $A \in \Gamma^{\infty}$,
$\mathrm{ZFC} \vdash_{\Omega} \mbox{\rm `` All
$\underTilde{\Sigma}^1_1(A)$-open games
are determined''}$.

Then for each $A\in \Gamma^{\infty}$, for each
$\Sigma^2_1$-formula $\phi(x)$; either
$\mathrm{ZFC} + \mathrm{CH} \vdash_{\Omega}\phi[A]$
or
$\mathrm{ZFC} + \mathrm{CH} \vdash_{\Omega}(\neg\phi)[A]$.\qed}
\medskip

Using the theorem on the determinacy of $\Gamma^{\infty}$-clopen
games one obtains the  converse.

\medskip
{\bf Theorem 2.4} {\em Suppose that there exists a proper class
of Woodin cardinals which are limits of Woodin cardinals. Let
$\Gamma^{\infty}$ be the set of all $A \subseteq \mathbb{R}$
such that $A$ is universally Baire.
Then the following are equivalent.
\begin{enumerate}

\item[\rm(1)]For each $A \in \Gamma^{\infty}$,
$\mathrm{ZFC} \vdash_{\Omega} \mbox{\rm `` All
$\underTilde{\Sigma}^1_1(A)$-open games
are determined''}$.
\item[\rm(2)] For each $A\in \Gamma^{\infty}$, for each
$\Sigma^2_1$-formula $\phi(x)$,  either
$\mathrm{ZFC} + \mathrm{CH} \vdash_{\Omega}\phi[A]$
or
$\mathrm{ZFC} + \mathrm{CH} \vdash_{\Omega}(\neg\phi)[A]$.\qed
\end{enumerate}}
\medskip

A set $A \subseteq \mathbb{N}$ is {\em $\Omega$-recursive}
if there exists a formula $\phi(x)$ such that for all
$k \in \mathbb{N}$, either
$\mathrm{ZFC}\vdash_{\Omega} \phi[k]
$ or $\mathrm{ZFC} \vdash_{\Omega} (\neg\phi)[k]$; and such that
$$
A = \{k \in \mathbb{N}{\large\;\mid\;} \mathrm{ZFC}\vdash_{\Omega} \phi[k]\}.
$$

The question of whether there exists a sentence $\Psi$ such that
for each $\Sigma^2_2$ sentence $\phi$, either
$\mathrm{ZFC} + \mathrm{CH} + \Psi
\vdash_{\Omega} \phi$, or $\mathrm{ZFC} + \mathrm{CH} + \Psi
\vdash_{\Omega} (\neg\phi)$,
and such that $\mathrm{ZFC} + \mathrm{CH} + \Psi$ is $\Omega$-consistent;
can be reformulated as:
\begin{quote}{\em Suppose there exists a
proper class of Woodin cardinals and that
$\mathrm{CH}$ holds.  Let $T$ be
the set of all $\Sigma^2_2$-sentences, $\phi$, such that
$$
V \vDash \phi.
$$
Can $T$ be $\Omega$-recursive?}
\end{quote}

\medskip
{\bf Theorem 2.5} {\em Suppose that there exists a proper class
of inaccessible limits of Woodin cardinals. Let
$\Gamma^{\infty}$ be the set of all $A \subseteq \mathbb{R}$
such that $A$ is universally Baire and suppose that
all $\Gamma^{\infty}$-open games
are determined.
Let $T_{\rm max}$ be the set of all $\Sigma^2_2$ sentences $\phi$
such that
$\mathrm{ZFC} + \mathrm{CH} + \phi$
is $\Omega$-consistent.

Then
$\mathrm{ZFC} + \mathrm{CH} + T_{\rm max}$
is $\Omega$-consistent.\qed}
\medskip

The following conjecture can be proved from rather
technical assumptions on the exsitence of an inner model theory for
the large cardinal hypothesis:
{\em $\kappa$ is $\delta$ supercompact where $\delta > \kappa$ and
$\delta$ is a Woodin cardinal.} The conjecture is:
\begin{quote}{\em
 Suppose that there exists a proper class
of supercompact  cardinals.
Let $T_{\rm max}$ be the set of all $\Sigma^2_2$ sentences $\phi$
such that
$\mathrm{ZFC} + \mathrm{CH} + \phi$
is $\Omega$-consistent.
Then $T_{\rm max}$ is $\Omega$-recursive.}
\end{quote}

While the plausibility of this conjecture is some evidence that $\Sigma^2_2$ absoluteness is possible, it
does not connect $\Sigma^2_2$ absoluteness with any determinacy hypothesis.

From inner model theory considerations any such determinacy hypothesis
must be beyond the
reach of superstrong cardinals.
In fact, Itay Neeman has defined a family of games whose (provable) determinacy
is arguably beyond the reach of superstrong cardinals; \cite{neeman}.

\section{Neeman games} \label{section 3}\setzero
\vskip-5mm \hspace{5mm }

For each formula $\phi(x_1,\ldots,x_n)$, let $X_{\phi}$ be the set
of all $a \in \{0,1\}^{\omega_1}$ such that there exists a closed,
unbounded set $C \subseteq \omega_1$ such that for all $\alpha_1 <
\cdots < \alpha_n$ in $C$,
$$
\<H(\omega_1),a,\in\> \vDash \phi[\alpha_1,\ldots, \alpha_n].
$$

The game given by $X_{\phi}$ is a {\em Neeman game}.
Are Neeman games determined?

Surprisingly the consistency strength of the determinacy of all
Neeman games is relatively weak.

\medskip
{\bf Lemma 3.1} {\em Suppose that all $\Delta^1_3$-clopen games are determined.
Then there exists $A \subseteq \omega_1$ such that in $L[A]$ if
$X \subseteq \{0,1\}^{\omega_1}$ is definable an $\omega$-sequence
of ordinals,
then $X$ is determined.\qed}
\medskip

One can easily introduce additional predicates for sets of reals.

For each formula, $\phi(x_1,\ldots,x_n)$, and for each
set $A \subseteq \mathbb{R}$ let $X_{(\phi,A)}$ be
the set of all $a \in \{0,1\}^{\omega_1}$ such that
there exists a closed, unbounded set $C \subseteq \omega_1$
such that for all $\alpha_1 < \cdots < \alpha_n$ in $C$,
$\<H(\omega_1),a,A,\in\> \vDash \phi[\alpha_1,\ldots, \alpha_n]$.
The game given by $X_{(\phi,A)}$ is an $A$-{\em Neeman} game.

\medskip
{\bf Definition 3.2} {\it $\diamond_{\mbox{\rm\tiny G}}$: For each
$\Sigma^2_2$ sentence, $\phi$, $V\vDash \phi$ if and only if
$V^{\mathrm{Coll}(\omega_1,\mathbb{R})} \vDash \phi$}.\qed
\medskip

The principle, $\diamond_{\mbox{\rm\tiny G}}$,  is a {\em generic}
form of $\diamond$. The next theorem gives a connection
between versions of $\Sigma^2_2$ absoluteness and determinancy
specifically the determinacy of
Neeman games. In this theorem it is the principle,
$\diamond_{\mbox{\rm\tiny G}}$, which plays the role
of $\mathrm{CH}$ in the theorem on $\Sigma^2_1$ absoluteness.

\medskip
{\bf Theorem 3.3} {\em Suppose that there exists a proper class
of supercompact  cardinals. Let
$\Gamma^{\infty}$ be the set of all $A \subseteq \mathbb{R}$
such that $A$ is universally Baire.
Then the following are equivalent.
\begin{enumerate}
\item[\rm(1)] For each $A \in \Gamma^{\infty}$,
$\mathrm{ZFC} + \diamond_{\mbox{\rm\tiny G}}\vdash_{\Omega}
\mbox{\rm `` All
$A$-Neeman games
are determined''}$.
\item[\rm(2)]  For each $A\in \Gamma^{\infty}$, for each
$\Sigma^2_2$-formula $\phi(x)$, either
$\mathrm{ZFC} + \diamond_{\mbox{\rm\tiny G}} \vdash_{\Omega}\phi[A]$
or
$\mathrm{ZFC} + \diamond_{\mbox{\rm\tiny G}} \vdash_{\Omega}(\neg\phi)[A]$.\qed
\end{enumerate}}
\medskip

We note the following trivial lemma which simply connects
the results here with the earlier ``evidence''
that $\Sigma^2_2$ absoluteness is possible; \mbox{cf.} the
discussion after Theorem~2.5.

\medskip

{\bf Lemma 3.4} {\em Suppose that there exists a proper class
of inaccessible limits of Woodin cardinals and
suppose that  for each
$\Sigma^2_2$-sentence $\phi$, either
$\mathrm{ZFC} + \diamond_{\mbox{\rm\tiny G}} \vdash_{\Omega}\phi$
or
$\mathrm{ZFC} + \diamond_{\mbox{\rm\tiny G}} \vdash_{\Omega}(\neg\phi)$.
Then for each $\Sigma^2_2$ sentence $\phi$ the following
are equivalent:
\begin{enumerate}
\item[\rm(1)] $\mathrm{ZFC} + \diamond_{\mbox{\rm\tiny G}} \vdash_{\Omega}\phi$;
\item[\rm(2)] $\mathrm{ZFC} + \mathrm{CH} +\phi$ is $\Omega$-consistent.
\qed
\end{enumerate}}

\medskip

The next theorem suggests that $\Sigma^2_2$ absoluteness conditioned
simply on $\diamond$ might actually follow from some large
cardinal hypothesis. Such a theorem would certainly be a striking
generalization of Theorem~1.1 and its proof might well
yield fundamental new insights into the combinatorics of
subsets of $\omega_1$.

\medskip
{\bf Theorem 3.5} {\em Suppose that there exists a proper class
of supercompact  cardinals. Let
$\Gamma^{\infty}$ be the set of all $A \subseteq \mathbb{R}$
such that $A$ is universally Baire and suppose that for each $A \in \Gamma^{\infty}$,
$\mathrm{ZFC}  \vdash_{\Omega} \mbox{\rm `` All
$A$-Neeman games
are determined''}$.
Then for each $A\in \Gamma^{\infty}$, for each
$\Sigma^2_2$-formula $\phi(x)$, either
$\mathrm{ZFC} + \diamond \vdash_{\Omega}\phi[A]$
or
$\mathrm{ZFC} + \diamond \vdash_{\Omega}(\neg\phi)[A]$.\qed}
\medskip

Given Theorem~3.5, the natural conjecture is that
Theorem~3.3 holds with $\diamond_{\mbox{\rm\tiny G}}$
replaced by $\diamond$. The missing ingredient in
proving such a conjecture seems to be a lack of
information on the nature of definable winning strategies
for Neeman games and more fundamentally on the lack of
any genuine determinacy proofs whatsoever for Neeman games.

In an exploration of the combinatorial aspects of
Neeman games it is useful to consider a wider class of games.
This class we now define.

For each formula, $\phi(x_1,\ldots,x_n)$, and for each
stationary set $S \subseteq \omega_1$ let $Y_{\phi}$ be
the set of all $a \in \{0,1\}^{\omega_1}$ such that
there exists a stationary set $S \subseteq \omega_1$
such that for all $\alpha_1 < \cdots < \alpha_n$ in $S$,
$\<H(\omega_1),a,\in\> \vDash \phi[\alpha_1,\ldots, \alpha_n]$.
The game given by $Y_{\phi}$ is a {\em stationary Neeman game}.

Can some large cardinal hypothesis imply that all stationary Neeman
games are determined?
Given the impossibility of $\Sigma^2_2({\cal I}_{_{\mathrm{NS}}})$-absoluteness, modulo
failure of the $\Omega$~Conjecture one would naturally conjecture
that the answer is ``no''. This is simply because there is
no apparent candidate for an absoluteness theorem which would
correspond to the (provable) determinacy of all stationary Neeman
games.

We define two games of length $\omega_1$. The first is a Neeman game
and the second is a stationary Neeman game.
Rather than have the moves be from $\{0,1\}$ it is more
convenient to have the moves be from $H(\omega_1)$.

\medskip

\noindent{\bf The canonical function game:} Player~I plays
$<a_{\alpha}:\alpha < \omega_1\>$
and Player~II plays
$<b_{\alpha}:\alpha < \omega_1\>$
subject to the rules: $a_{\alpha +1} \subset \alpha\times\alpha$
and $b_{\alpha}$ is a countable ordinal.

Player~I wins if there exists a set $A \subset \omega_1\times\omega_1$
such that
$A$ is a wellordering of $\omega_1$
and such that there exists a closed unbounded set $C \subset \omega_1$
such that for all $\alpha \in C$:
$a_{\alpha+1} = A\cap (\alpha\times \alpha)$ and
$b_{\alpha} < \mbox{\rm rank}(a_{\alpha+1})$.

\medskip

\noindent{\bf The stationary canonical function game:} Player~I plays
$<a_{\alpha}:\alpha < \omega_1\>$
and Player~II plays
$<b_{\alpha}:\alpha < \omega_1\>$
subject to the rules: $a_{\alpha +1} \subset \alpha\times\alpha$
and $b_{\alpha}$ is a countable ordinal.

Player~I wins if there exists a set $A \subset \omega_1\times\omega_1$
such that
$A$ is a wellordering of $\omega_1$
and such that there exists a statationary set $S \subset \omega_1$
such that for all $\alpha \in S$:
$a_{\alpha+1} = A\cap (\alpha\times \alpha)$ and
$b_{\alpha} < \mbox{\rm rank}(a_{\alpha+1})$.
\medskip

In models where $L$-like condensation principles hold
these games are easily seen to be determined.

\medskip
{\bf Lemma 3.6} {\em Suppose $\diamond$ holds. Then Player~II has a
winning strategy in the canonical function game.\qed}
\medskip

{\bf Lemma 3.7} {\em Suppose $\diamond^+$ holds. Then Player~II has
a winning strategy in the stationary canonical function game.\qed}
\medskip

In contrast to the previous lemma, the following theorem shows
that it is consistent that Player~I has a winning strategy in
the stationary canonical function game, at least if fairly
strong large cardinal hypotheses are assumed to be consistent.

\medskip
{\bf Theorem 3.8}{\em
\hspace{1mm} Suppose there is a huge cardinal. Then there
is a partial order, $\mathbb{P}$, such that in $V^{\mathbb{P}}$,
Player~I has a winning strategy
in the  stationary canonical function game.\qed}
\medskip

These two results strongly suggest that no large cardinal hypothesis
can imply that the stationary canonical function game is determined.
In fact from consistency of a relatively weak large cardinal
hypothesis, one does obtain the consistency that the stationary
canonical function game is not determined. Note that if the
stationary canonical function game is not determined then
every function,
$f:\omega_1 \to \omega_1$,
is bounded by a canonical function on a stationary set and
so the consistency of {\em some} large cardinal hypothesis is necessary.

\medskip
{\bf Theorem 3.9} {\em \hspace{1mm}
Suppose there is a measurable cardinal.
Then there
is a partial order, $\mathbb{P}$, such that in $V^{\mathbb{P}}$
the stationary canonical function game  is not determined.\qed}
\medskip

There are many open problems about the
canonical function games. Here are several.
\begin{enumerate}

\item {\em Is it consistent that Player~I has a winning strategy
in the canonical
function game?  }

\item {\em Is it consistent that Player~II does not
have a winning strategy in the canonical function game? }

\item {\em Is it consistent that Player~I has a winning strategy
in the stationary canonical
function game on each stationary set?  }

\item {\em How strong is the assertion that Player~I has a winning
strategy in the stationary canonical function game?}

\end{enumerate}

For each formula, $\phi(x_1,\ldots,x_n)$,
for each sequence
$$
{\cal S} = \<S_{\alpha}:\alpha < \omega_1\>
$$
of pairwise disjoint stationary subsets of $\omega_1$ and
that $A \subseteq \mathbb{R}$, let $Y^{\cal S}_{(\phi,A)}$
be the set of all $a \in \{0,1\}^{\omega_1}$ such that
there exists a stationary set $S \subseteq \omega_1$
such that for all $\alpha_1 < \cdots < \alpha_n$ in $S$,
$$
\<H(\omega_1),a,A,\in\> \vDash \phi[\alpha_1,\ldots, \alpha_n],
$$
 and such that $S\cap S_{\alpha}$ is stationary for all $\alpha < \omega_1$.

\medskip
{\bf Theorem 3.10} {\em Suppose that there exists a proper class
of supercompact  cardinals. Let
$\Gamma^{\infty}$ be the set of all $A \subseteq \mathbb{R}$
such that $A$ is universally Baire.

Suppose that $A \in \Gamma^{\infty}$, $\phi(x_1,\ldots,x_n)$
is a formula and that
$$
\mathrm{ZFC}  \vdash_{\Omega} \mbox{\rm `` The Neeman game
$X_{(\phi,A)}$ is determined".}
$$

Then either:
\begin{enumerate}
\item[\rm(1)] $\mathrm{ZFC}  \vdash_{\Omega}
\mbox{\rm `` I wins the game $X_{(\phi,A)}$''}$, or;
\item[\rm(2)]  $\mathrm{ZFC}  \vdash_{\Omega}
\mbox{\rm `` For all ${\cal S}$, II wins the game $Y^{\cal S}_{(\phi,A)}$''}$.
\qed
\end{enumerate}}
\medskip

The determinacy hypothesis:
{\em
All Neeman games are determined};
is relatively weak in consistency strength (the consistency strength is
at most that of the existence of a Woodin cardinal which is a limit
of Woodin cardinals).
However
the determinacy hypothesis:
\begin{quote}{\em
For each formula $\phi$, either Player~I wins the game $X_{\phi}$,
or for each sequence,
$$
{\cal S} = \<S_{\alpha}:\alpha < \omega_1\>,
$$
of pairwise disjoint stationary subsets of $\omega_1$,
Player~II wins $Y^{\cal S}_{(\phi,\emptyset)}$};
\end{quote}
seems plausibly very strong.

\label{lastpage}

\end{document}